\documentclass[12pt]{amsart}
\usepackage{amssymb}
\newcommand{\R}{\mathbb{R}}
\newcommand{\N}{\mathbb{N}}
\newcommand{\Z}{\mathbb{Z}}
\newcommand{\C}{\mathbb{C}}

\newcommand{\F}{\mathcal{F}}

\def\endproof{\quad \hfill$\blacksquare$\vspace{0.15cm}\\}
\newcommand{\ds}{\displaystyle}

\newtheorem{propo}{ Proposition}
\newtheorem{lemme}{Lemma}

\newtheorem{theorem}{Theorem}

\begin{document}

\title [Uncertainty principles for the Fourier transforms in Quantum Calculus]
 { Uncertainty principles for the Fourier transforms in Quantum Calculus}
\author{N\'eji Bettaibi \quad \&\quad Ahmed Fitouhi \quad \&\quad Wafa Binous  }
\address{N. Bettaibi. Institut
Pr\'eparatoire  aux \'Etudes d'Ing\'enieur de Nabeul, 8000 Nabeul,
Tunisia.}
 \email{Neji.Bettaibi@ipein.rnu.tn}

\address{A. Fitouhi. Facult\'e des Sciences de Tunis,
 1060 Tunis, Tunisia.}
\email{Ahmed.Fitouhi@fst.rnu.tn}

\begin{abstract} Some properties of the $q$-Fourier-sine  transform are studied and
 $q$-analogues of the Heisenberg uncertainty principle is
derived for the $q$-Fourier-cosine  transform studied in \cite{FB}
and for the $q$-Fourier-sine  transform.

\end{abstract}
 \maketitle
\section{Introduction}

One of the basic  principles in classical Fourier analysis is the
impossibility to find a function $f$ being arbitrarily well
localized together with its Fourier transform $\widehat{f}$. There
are many ways to get this statement precise. The most famous of
them is the so called Heisenberg uncertainty principle, a
consequence of Cauchy-Schwarz's inequality which states that for
$f\in L^2(\R)$,
\begin {equation}
\left(\int_{-\infty}^\infty x^2\mid f(x)\mid^2dx\right)
\left(\int_{-\infty}^\infty \lambda^2\mid
\widehat{f}(\lambda)\mid^2d\lambda\right) \geq
\frac{1}{4}\left(\int_{-\infty}^\infty \mid f(x)\mid^2dx \right)^2
\end{equation}
with equality only if $f(x)$ is almost everywhere equal to a
constant multiple of $e^{-px^2}$ for some $p>0$. Here
$$\widehat{f}(\lambda)=\frac{1}{\sqrt{2\pi}}\int_{-\infty}^\infty f(x)e^{-ix\lambda}dx.$$
In this paper we shall prove that similarly to the classical
theory, a nonzero function and its $q$-Fourier ($q$-Fourier-cosine
and $q$-Fourier-sine) transform cannot both be sharply localized.
For this purpose we will prove a $q$-analogue of the Heisenberg
uncertainly principle. This paper is organized as follows: in
Section 2, we present some preliminaries results and notations
that will be  useful in the sequel. In Section 3, we study some
$q$-harmonic  results and state  $q$-analogues of the Heisenberg
uncertainly principle.

\section{notations and preliminaries}
\indent Throughout this paper, we will fix $q \in ]0,1[$ such that
$\ds \frac{Log(1-q)}{Log(q)}\in\Z$ .
 We recall some usual notions and notations used in the
$q$-theory (see \cite {GR} and \cite{KC}). We refer to the book by
G. Gasper and M. Rahmen \cite{GR}, for the definitions, notations
and properties of   the $q$-shifted factorials and the
$q$-hypergeometric functions.\\

 We note\\
 $\ds \R_q=\{\pm q^n~~:~~n\in\Z\}$ and  $\ds \R_{q,+}=\{
 q^n~~:~~n\in\Z\}$.\\
 We also denote
 \begin {equation}
 [x]_q={{1-q^x}\over{1-q}},~~~~~ x\in \C
\end {equation}
and
\begin {equation}
    [n]_q! ={{(q;q)_n}\over
 {(1-q)^n}},  ~~~~~~ n\in \N.
\end {equation}
 The $q$-derivatives
$D_qf$ and $D_q^+f$ of a function $f$ are given by
\begin {equation}\label {21}
(D_qf)(x)={{f(x)-f(qx)}\over{(1-q)x}},~~
(D_q^+f)(x)={{f({q^{-1}x})-f(x)}\over{(1-q)x}}, ~~ { if}~~
x\not=0,
\end {equation}
$(D_qf)(0)=f'(0)$ and $(D_q^+f)(0)=q^{-1}f'(0)$  provided $f'(0)$
 exists.\\
 The $q$-Jackson integrals from $0$ to $a$ and from $0$ to $\infty$ are
defined by (see \cite {Jac})
\begin {equation} \label {26}
\int_0^{a}{f(x)d_qx} =(1-q)a\sum_{n=0}^{\infty}{f(aq^n)q^n},
\end {equation}
\begin {equation} \label{int}
\int_0^{\infty}{f(x)d_qx}
=(1-q)\sum_{n=-\infty}^{\infty}{f(q^n)q^n},
\end {equation}
provided the sums converge absolutely.\\The $q$-Jackson integral
in a generic interval $[a,b]$ is given  by (see \cite {Jac})
\begin {equation}
\int_a^{b}{f(x)d_qx} =\int_0^{b}{f(x)d_qx}-\int_0^{a}{f(x)d_qx}.
\end {equation}
 The improper integral is defined in the following way
(see \cite {KoC})
\begin {equation}
\int_0^{{\infty}\over A}{f(x)d_qx}
  =(1-q)\sum_{n=-\infty}^{\infty}{f\left({{q^n}\over A}\right){{q^n}\over A}
  }.
\end {equation}
 We remark that for $n\in \Z$, we have
\begin {equation}\label {30}
\int_0^{{\infty}\over {q^n}}{f(x)d_qx}=\int_0^{\infty}{f(x)d_qx}.
\end {equation}
 The  $q$-integration by parts is given for suitable functions $f$ and $g$ by
\begin {equation}
 \int_a^{b}{g(x)D_qf(x)d_qx}=f(b)g(b)-f(a)g(a)-\int_a^{b}{f(qx)D_qg(x)d_qx}.
\end {equation}
{\bf Remark} A second $q$-analogue of the integration by parts
theorem is given for a suitable function $f$ and $g$ by (see
\cite{KoT})
\begin {equation}\label{ipp}
 \int_a^{b}{g(x)D_qf(x)d_qx}=f(b)g(q^{-1}b)-f(a)g(q^{-1}a)-\int_a^{b}{f(x)D_q^+g(x)d_qx}.
\end {equation}
\begin {propo}The $q$-analogue of the integration theorem by change of variable
 for  $u(x)=\alpha x^{\beta}, \quad \alpha \in \C$ and
$\beta >0$ is as follows
\begin {equation}
\int_{u(a)}^{u(b)}{f(u)d_qu} =
\int_a^b{f(u(x))}D_{q^{\frac{1}{\beta}}}u(x)d_{{q^{\frac{1}{\beta}}}}x.
\end {equation}
\end {propo}

\indent Jackson \cite {Jac} defined the $q$-analogue of the Gamma
function by
\begin {equation}
\Gamma_q (x) ={(q;q)_{\infty}\over{(q^x;q)_{\infty}}}(1-q)^{1-x} ,
\qquad x\not={0},{-1},{-2},\ldots.
\end {equation}
It is well known that it satisfies
\begin {equation}
\Gamma_q(x+1)=[x]_q \Gamma_q(x),\quad \Gamma_q(1)=1~~~~{\rm
and}~~~~ \lim_{q\rightarrow1^-}\Gamma_q(x)=\Gamma(x),~~ \Re(x)>0.
\end {equation}

 The third $q$-Bessel function (see \cite{Ism, KS}) is given and denoted by M.
E. H. Ismail  as
\begin {equation}
J_{\alpha}(z;q^2)=\frac{z^{\alpha}}{(1-q^2)^{\alpha}\Gamma_{q^2}(\alpha+1)}
        {_1\varphi_1(0;q^{2\alpha + 2};q^2,q^2z^2)}.
\end  {equation}
It verifies for $\alpha,~~ \beta$ reals (see \cite{KS})
\begin {equation}
J_{\alpha}(q^{\beta};q^2) =J_{\beta}(q^{\alpha};q^2),
\end {equation}
and we have the following orthogonality relation (see \cite{KS})
\begin {equation}\label{ort}
\sum_{k=-\infty}^\infty
q^{2k}q^{n+m}J_{n+k}(x;q^2)J_{m+k}(x;q^2)=\delta_{n,m},~~\mid
x\mid <q^{-1},~~~~n,~m\in\Z.
\end {equation}
Moreover, if $\ds \alpha >-1$, we have (see \cite{FKB}, \cite{KS}
),
 \begin {equation} \label{asj}
    \forall x\in \R_{q,+},~~~~\mid
J_{\alpha}(x;q^2)\mid\leq\frac{(-q^2;q^2)_{\infty}(-q^{2(\alpha+1)};q^2)_{\infty}}{(q^2;q^2)_{\infty}}\left\{
\begin{array}{ccc}
  1 & {      \textrm{if}} & x\leq1 \\
  q^{(\frac{Logx}{Logq})^2} & {      \textrm{if}} & x>1.
\end{array} \right.
\end{equation}
and
\begin {equation} \forall \nu\in \R,
~~J_{\alpha}(x;q^2)=o(x^{-\nu})~~~~{ \textrm{as}}
~~~~x\rightarrow+\infty ~~~~~{      \textrm{in}}~~
\R_{q,+}.\end{equation}

The $q$-trigonometric functions $q$-cosine and $q$-sine are
defined by ( see \cite {FB}):
\begin {equation}
 \cos(x;q^2)=~~
_1\varphi_1\biggl(0,q;q^2,(1-q)^2x^2\biggr)=\sum_{n=0}^{\infty}(-1)^nq^{n(n-1)}
{{x^{2n}}\over{[2n]_q!}}
\end {equation}
and
\begin {equation}
 \sin(x;q^2)=x~~_1\varphi_1\biggl(0,q^3;q^2,(1-q)^2x^2\biggr)=\sum_{n=0}^{\infty}(-1)^nq^{n(n-1)}
{{x^{2n+1}}\over{[2n+1]_q!}}.
\end {equation}

Note that we have the relations
\begin {equation} \label{cos}
\cos(x;q^2) =
\frac{\Gamma_{q^2}(\frac{1}{2})}{q(1+q^{-1})^{\frac{1}{2}}}x^{\frac{1}{2}}J_{-\frac{1}{2}}(\frac{1-q}{q}x;q^2),
\end {equation}

\begin {equation} \label {sj}
\sin(x;q^2) =
\frac{\Gamma_{q^2}(\frac{1}{2})}{(1+q^{-1})^{\frac{1}{2}}}x^{\frac{1}{2}}J_{\frac{1}{2}}(\frac{1-q}{q}x;q^2)
\end {equation}
and they verify

\begin {equation}
D_q\cos(x;q^2)=-\frac{1}{q}\sin(qx;q^2)
\end {equation}
and
\begin {equation}
D_q\sin(x;q^2)=\cos(x;q^2).
\end {equation}
\section{$q$-Uncertainly principle}
We define the $q$-Fourier-cosine and the $q$-Fourier-sine as ( see
\cite{FB} and \cite{KS})
\begin {equation}
\F_q(f)(x)=c_q\int_0^\infty f(t)\cos(xt;q^2)d_qt
\end {equation}
and
\begin {equation}
_q\F(f)(x)=c_q\int_0^\infty f(t)\sin(xt;q^2)d_qt,
\end {equation}
where
\begin {equation}
c_q= \frac{(1+q^{-1})^{\frac{1}{2}}}{\Gamma_{q^2}(\frac{1}{2})}.
\end {equation}
It was shown in \cite{D} that $\F_q$ is an isomorphism of
$L_q^2(\R_{q,+})$ and we have $\F_q^{-1}= \F_q$ and the following
Plancherel formula: $$ \parallel \F_q(f)\parallel_{q,2}=\parallel
f\parallel_{q,2},~~~~f\in L_q^2(\R_{q,+}),$$ where
$L_q^n(\R_{q,+})$ is the set of functions defined on $\R_{q,+}$
such that $\ds \int_0^\infty \mid f(t)\mid^nd_qt < \infty$,
equipped with the norm $\ds \parallel f\parallel_{q,n}=\left(
\int_0^\infty \mid f(t)\mid^nd_qt\right)^{\frac{1}{n}}$. \\
The $q$-Fourier-sine verifies the following properties.
\begin{propo}
For $f\in L_q^1(\R_{q,+})$, we have\\
1) $\ds \forall \lambda \in \R_{q,+},
\mid~_q\F(f)(\lambda)\mid\leq
\frac{(1+q^{-1})^{1/2}}{\Gamma_{q^2}(1/2)
(q;q)_\infty^2}.\parallel f\parallel_{q,1};$\\
2) $\ds \lim_{\lambda\rightarrow \infty}~_q\F(f)(\lambda)=0.$
\end{propo}
\proof Using the inequality (see \cite{FB})
$$ \mid \sin(x;q^2)\mid \leq \frac{1}{(q;q)_\infty^2},~~x\in \R_q,$$
we obtain
$$ \mid f(t)\sin(\lambda t;q^2)\mid \leq \frac{1}{(q;q)_\infty^2}\mid f(t)\mid ,~~\lambda, t\in \R_q.$$
Which gives, after integration, the first inequality and together
with the Lebesgue theorem it gives the limit. \endproof \\
In the following proposition, we shall try to prove a Plancherel
formula for the $q$-Fourier-sine transform. We begin by the
following useful result:
\begin {lemme}
 For all
$x, y\in \R_{q,+}$, we have
\begin {equation}
\sqrt{xy}\int_0^\infty \sin(xt;q^2)\sin(yt;q^2)d_qt=
\frac{q^2\Gamma_{q^2}^2(\frac{1}{2})}{(1+q^{-1})(1-q)}\delta_{x,y}.
\end {equation}
\end {lemme}
\proof Let $x=q^n$ and $y=q^m$, $m,n \in\Z$ be two elements of
$\R_{q,+}$. The orthogonality relation (\ref{ort}) leads to
$$\sum_{k=-\infty}^\infty
q^{2k}q^{n+m}J_{n+k}(q^{\frac{1}{2}};q^2)J_{m+k}(q^{\frac{1}{2}};q^2)=\delta_{n,m},$$
which is equivalent to
$$\sum_{k=-\infty}^\infty
q^{2k}q^{n+m}J_{\frac{1}{2}}(q^{n+k};q^2)J_{\frac{1}{2}}(q^{m+k};q^2)=\delta_{n,m}.$$
Using the relation (\ref{sj}), we obtain
$$\sum_{k=-\infty}^\infty \frac{(1-q)(1+q^{-1})}{q\Gamma_{q^2}^2(\frac{1}{2})}
q^{\frac{n+m}{2}}q^{k}\sin\left(\frac{q}{1-q}q^{n+k};q^2\right)\sin\left(\frac{q}{1-q}q^{m+k};q^2\right)=\delta_{n,m}.$$
Then $$\frac{1+q^{-1}}{q\Gamma_{q^2}^2(\frac{1}{2})}
q^{\frac{n+m}{2}} \int_0^\infty \sin
\left(\frac{q}{1-q}q^{n}t;q^2\right)\sin\left(\frac{q}{1-q}q^{m}t;q^2\right)d_qt=\delta_{n,m}.$$
The change of variable $u= \frac{q}{1-q}t$ gives
$$\frac{(1+q^{-1})(1-q)}{q^2\Gamma_{q^2}^2(\frac{1}{2})}
q^{\frac{n+m}{2}} \int_0^\infty \sin
\left(q^{n}t;q^2\right)\sin\left(q^{m}t;q^2\right)d_qt=\delta_{n,m}.$$
Thus
$$\sqrt{xy}\int_0^\infty \sin
\left(xt;q^2\right)\sin\left(yt;q^2\right)d_qt=\frac{q^2\Gamma_{q^2}^2(\frac{1}{2})}{(1+q^{-1})(1-q)}\delta_{x,y}.$$
\endproof
\begin{propo}\label{spl} 1) For
$f\in L_q^2(\R_{q,+})$, we have $_q\F(f)\in L_q^2(\R_{q,+})$ and
$$\ds
\parallel _q\F(f)\parallel_{q,2}= q\parallel
f\parallel_{q,2}.$$
 2)$_q\F$ is an isomorphism of
$L_q^2(\R_{q,+})$ and  $\ds (_q\F)^{-1}=\frac{1}{q^2}~~_q\F$.
\end{propo}
\proof 1) For $x\in \R_{q,+}$, we have
\begin{eqnarray*}
_q\F(f)(x)&=&c_q\int_0^\infty f(t)\sin(xt;q^2)d_qt\\
     &=& c_q (1-q)\sum_{n=-\infty}^\infty q^nf(q^n)\sin(xq^n;q^2).
\end{eqnarray*}
So, for $x\in \R_{q,+}$,
\begin{eqnarray*}
\left( _q\F(f)(x)\right)^2&=& c_q^2 (1-q)^2\sum_{n=-\infty}^\infty
q^{2n}f^2(q^n)\sin^2(xq^n;q^2)\\
&+& c_q^2 (1-q)^2\sum_{n,m=-\infty, n\neq m}^\infty
q^{m+n}f(q^m)f(q^n)\sin(xq^m;q^2)\sin(xq^n;q^2) .
\end{eqnarray*}
By integration, we obtain
\begin{eqnarray*}
\int_0^\infty \left( _q\F(f)(x)\right)^2d_qx &=& c_q^2
(1-q)^2\int_0^\infty \sum_{n=-\infty}^\infty
q^{2n}f^2(q^n)\sin^2(xq^n;q^2)d_qx\\
&+& c_q^2 (1-q)^2\int_0^\infty \sum_{n,m=-\infty, n\neq m}^\infty
q^{m+n}f(q^m)f(q^n)\sin(xq^m;q^2)\sin(xq^n;q^2)d_qx.
\end{eqnarray*}
The previous lemma, the relation (\ref{asj}) and Fubini's theorem
imply that we can exchange  the integral and the sum signs and we
have:
\begin{eqnarray*}
\int_0^\infty \left( _q\F(f)(x)\right)^2d_qx &=& c_q^2 (1-q)^2
\sum_{n=-\infty}^\infty
q^{n}f^2(q^n) q^n\int_0^\infty\sin^2(xq^n;q^2)d_qx\\
&=& c_q^2
(1-q)^2\frac{q^2\Gamma_{q^2}^2(\frac{1}{2})}{(1+q^{-1})(1-q)}
\sum_{n=-\infty}^\infty
q^{n}f^2(q^n)\\
&=&q^2 \int_0^\infty f^2(t)d_qt.
\end{eqnarray*}
2) Using the same arguments, we can see that for $y\in\R_{q,+}$,
we have
\begin{eqnarray*}
\int_0^\infty~ _q\F(f)(x)\sin(xy;q^2)d_qx &=&
c_q(1-q)\sum_{n=-\infty}^\infty q^nf(q^n)\int_0^\infty
\sin(xq^n;q^2)\sin(xy;q^2)d_qx\\
&=& \frac{q^2}{c_q}  f(y).
\end{eqnarray*}
\endproof
The following result gives a relation between the
$q$-Fourier-cosine and the $q$-Fourier-sine.
\begin{lemme}\label{cfs}
For $f\in L_q^1(\R_{q,+})$ such that  $D_qf\in L_q^1(\R_{q,+})$,
we have:\\ 1)
\begin {equation}
_q\F(D_qf)(\lambda)=-\frac{\lambda}{q}
\F_q(f)\left(\frac{\lambda}{q}\right), ~~\lambda\in \R_{q,+}.
\end{equation}
2) Additionally, if $f(0)=0$  then
\begin {equation}
\F_q(D_qf)(\lambda)=\frac{\lambda}{q^2}~~ _q\F(f)(\lambda),
~~\lambda\in \R_{q,+}.
\end{equation}
\end{lemme}

\proof Since $f$ is in $L_q^1(\R_{q,+})$ then for all $\lambda\in
\R_{q,+}$, $f(t)\sin(\lambda t;q^2)$ and  $f(t)\sin(\lambda
t;q^2)$ tend to $0$ as $t$ tends to $\infty$. So by
$q$-integrations by parts, we obtain
\begin{eqnarray*}
 _q\F(D_qf)(\lambda)&=& c_q\int_0^\infty D_qf(t)\sin(\lambda t;q^2)d_qt\\
  &=& -c_q\int_0^\infty \lambda f(qt)\cos(\lambda t;q^2)d_qt\\
  &=& -\frac{\lambda}{q} \F_q(f)\left(\frac{\lambda}{q}\right)
\end{eqnarray*}
and
\begin{eqnarray*}
 \F_q(D_qf)(\lambda)&=& c_q\int_0^\infty D_qf(t)\cos(\lambda t;q^2)d_qt\\
  &=& \frac{c_q}{q}\int_0^\infty \lambda f(qt)\sin(q\lambda t;q^2)d_qt\\
  &=& \frac{\lambda}{q^2}~~ _q\F(f)(\lambda)
\end{eqnarray*}
\endproof
Now, we are in a situation to state a $q$-analogues of the
Heisinberg uncertainty principle.
\begin{theorem}
Let $f$ be in $L_q^2(\R_{q,+})$ such that $D_qf$ is in
$L_q^2(\R_{q,+})$. Then
\begin {equation}\label{upc}
\left(\int_0^\infty t^2\mid
f(t)\mid^2d_qt\right)^{1/2}\left(\int_0^\infty x^2\mid
\F_q(f)(x)\mid^2d_qx\right)^{1/2}\geq
\frac{q}{q^{\frac{3}{2}}+1}\parallel f
\parallel_{q,2}^2
\end {equation}
\end{theorem}
\proof First, using the previous lemma and Proposition \ref{spl},
we have
\begin{eqnarray*}
\int_0^\infty \lambda^2\mid \F_q(f)(\lambda)\mid^2d_q\lambda &=&
1/q\int_0^\infty  \mid
\frac{\lambda}{q}\F_q(f)(\frac{\lambda}{q})\mid^2d_q\lambda\\
&=& 1/q \int_0^\infty  \mid~~
_q\F(D_qf)(\lambda)\mid^2d_q\lambda\\
&=&q\int_0^\infty  \mid D_qf(t)\mid^2d_qt.
\end{eqnarray*}
The relation
$$D_q(f\overline{f})(t)= D_qf(t)\overline{f} +
f(qt)D_q\overline{f}(t)$$ leads to
\begin{eqnarray*}
 \left|\int_0^\infty tD_q(f\overline{f})(t)d_qt\right|&\leq&
 \int_0^\infty \left|tD_qf(t)\overline{f}(t)\right|d_qt + \int_0^\infty
 \left|tf(qt)D_q\overline{f}(t)  \right|d_qt\\
 &\leq& \left(\int_0^\infty
 \left|t\overline{f}(t)\right|^2d_qt\right)^{1/2}\left(\int_0^\infty
 \left|D_qf(t)\right|^2d_qt\right)^{1/2}\\ &+& \left(\int_0^\infty
 \left|tf(qt)\right|^2d_qt\right)^{1/2}\left(\int_0^\infty
 \left|D_q\overline{f}\right|^2d_qt\right)^{1/2}\\
 &=& \left(\frac{q^{3/2}+1}{q^2}\right) \left(\int_0^\infty t^2\left|f(t)\right|^2d_qt\right)^{1/2}\left(\int_0^\infty
 x^2\left|\F_q(f)(x)\right|^2d_qx\right)^{1/2}.
\end{eqnarray*}
On the other hand, since $f$ is in $L_q^2(\R_{q,+})$ then
$t|f(t)|^2$ tends to $0$ as $t$ tends to $\infty$ in $\R_{q,+}$.
So by $q$-integration by parts, we obtain
$$\int_0^\infty tD_q(f\overline{f})(t)d_qt = -\int_0^\infty
|f(qt)|^2d_qt=-1/q\int_0^\infty |f(t)|^2d_qt.$$

Finally
$$\left(\int_0^\infty t^2\mid
f(t)\mid^2d_qt\right)^{1/2}\left(\int_0^\infty x^2\mid
\F_q(f)(x)\mid^2d_qx\right)^{1/2}\geq
\frac{q}{q^{\frac{3}{2}}+1}\parallel f
\parallel_{q,2}^2.$$
\endproof
Similarly, we have an uncertainty principle for the
$q$-Fourier-sine transform.

\begin{theorem}Let $f$ be in $L_q^2(\R_{q,+})$ such that $D_qf$ is in
$L_q^2(\R_{q,+})$ and $f(0)=0$. Then
\begin {equation}\label{ups}
\left(\int_0^\infty t^2\mid
f(t)\mid^2d_qt\right)^{1/2}\left(\int_0^\infty \lambda^2\mid
_q\F(f)(\lambda)\mid^2d_q\lambda\right)^{1/2}\geq
\frac{q}{q^{\frac{-3}{2}}+1}\parallel f
\parallel_{q,2}^2.
\end {equation}
\end{theorem}
\proof Owing to Lemma \ref{cfs} and the Plancherel formula, we
have
\begin{eqnarray*}
\int_0^\infty \lambda^2 \mid ~ _q\F(f)(\lambda)\mid^2d_q\lambda
&=& q^4 \int_0^\infty  \mid
\F_q(D_qf)(\lambda)\mid^2d_q\lambda\\
&=&q^4\int_0^\infty  \mid D_qf(t)\mid^2d_qt.
\end{eqnarray*}

Using the same steps of the previous proof, we have

\begin{eqnarray*}
1/q\int_0^\infty |f(t)|^2d_qt
 &\leq&  \left(\int_0^\infty
\left|t\overline{f}(t)\right|^2d_qt\right)^{1/2}\left(\int_0^\infty
\left|D_qf(t)\right|^2d_qt\right)^{1/2}\\ &+& \left(\int_0^\infty
\left|tf(qt)\right|^2d_qt\right)^{1/2}\left(\int_0^\infty
\left|D_q\overline{f}\right|^2d_qt\right)^{1/2}\end{eqnarray*}
 $$=
\frac{1+q^{-3/2}}{q^2}\left(\int_0^\infty t^2\mid
f(t)\mid^2d_qt\right)^{1/2}\left(\int_0^\infty \lambda^2\mid
_q\F(f)(\lambda)\mid^2d_q\lambda\right)^{1/2}.$$

\endproof
{\bf Remark.} Note that when $q$ tends to $1$, the inequalities
(\ref{upc}) and (\ref{ups}) tend at least formally to the
corresponding classical ones.
\section{Uncertainty Principle in Hilbert space}
For $A$ and $B$ operators on a Hilbert space $H$, with domains
$D(A)$ and $D(B)$ respectively, we note\\
 $\ds [A,B]=AB-BA$  ~~ and ~~ $\ds [A,B]_q=qAB-BA$.
 The commutator $[A,B]$ and the $q-$ commutator $[A,B]_q$ are
 both defined on $D[A,B]=D(AB)\cap D(BA)$, where $\ds D(AB)=\{u\in D(B):Bu\in
 D(A)\}$ and likewise for $D(BA)$. Let us begin by the following
well-known result:
\begin{lemme}{\bf \rm(Cauchy-Schwarz's inequality)} For $x$ and
$y$ in the Hilbert space  $H$ the following inequality
\begin {equation}
|<x,y>|\leq \parallel x\parallel \parallel y\parallel
\end{equation}
holds.
\end{lemme}
Using this lemma, one can prove easily the following proposition,
which gives the uncertainty principle for normal operators.
\begin{propo} For $s\geq 0$, note $\ds [A,B]_s=sAB-BA$.
If $A$ and $B$ are  operators on the Hilbert space $H$, then for
all $u\in D[A,B]$, we have
\begin {equation}
\parallel Au\parallel \parallel B^*u\parallel+s\parallel
Bu\parallel \parallel A^*u\parallel\geq~~ <[A,B]_su,u>|.
\end {equation}
In addition, if $A$ and $B$ are normal on $H$, we obtain
\begin {equation}
\parallel Au\parallel \parallel Bu\parallel\geq\frac{1}{1+s}
|<[A,B]_s
u,u>|.
\end {equation}
\end{propo}
\proof For $s\geq0$, using the lemma  we have , and $u\in
D([A,B])$
\begin{eqnarray*}
|<(sAB-BA)u,u>|&=&|s<ABu,u>-<BAu,u>|\\&=&|s<Bu,A^*u>-<Au,B^*u>|\\
&\leq& |s<Bu,A^*u>| + |<Au,B^*u>|\\&\leq& s\parallel  Bu \parallel
\parallel A^*u\parallel + \parallel  Au \parallel
\parallel B^*u\parallel.
\end{eqnarray*}
Additionally, if  the operators $A$ and $B$ are normal, we obtain
$$|<(sAB-BA)u,u>|\leq(1+s) \parallel Au\parallel \parallel
Bu\parallel.$$
\endproof

Now, take $\ds H= \{f\in L_q^2(\R_{q,+}): f(0)=0\}$, $ A:
f(x)\mapsto xf(x)$, $ B: f(x)\mapsto iD_qf(x)$ and $ C:
f(x)\mapsto iD_q^+f(x)$. $A$ is self-adjoint on $H$,  $B^*= C$ (
according to (\ref{ipp}))  and
$$D[A,B]=D[A,C]=\{f\in H: xf, D_qf\in L_q^2(\R_{q,+}\}.$$
For all $f\in D[A,B]$, we have\\
$$[A,B]_qf(x)=-if(x)$$ and  $$[A,C]f(x)=-iq^{-1}f(q^{-1}x).$$
So, for all $f\in D[A,B]$, we have
\begin {eqnarray}\label{uph}
\parallel f\parallel_{q,2}^2&\leq& (\parallel D_q^+f\parallel_{q,2}+q \parallel
D_qf\parallel_{q,2}) \left(\int_0^\infty
x^2|f(x)|^2d_qx\right)^{1/2}\\
&=&\frac{1+q^{3/2}}{q^{1/2}}\parallel D_qf\parallel_{q,2}
\left(\int_0^\infty x^2|f(x)|^2d_qx\right)^{1/2}
\end {eqnarray}
and
\begin {eqnarray}
\left|\int_0^\infty f(x)\overline{f(qx)}d_qx\right|&\leq&
(\parallel D_q^+f\parallel_{q,2}+ \parallel D_qf\parallel_{q,2})
\left(\int_0^\infty
x^2|f(x)|^2d_qx\right)^{1/2}\\
&=&(1+\frac{1}{q^{1/2}})\parallel
D_qf\parallel_{q,2} \left(\int_0^\infty
x^2|f(x)|^2d_qx\right)^{1/2}.
\end {eqnarray}
{\bf Remark.} Using the fact that $\ds \parallel
D_qf\parallel_{q,2} = 1/\sqrt{q}\left(\int_0^\infty
\lambda^2|\F_q(f)(\lambda)|^2d_q\lambda\right)^{1/2}$, one can see
that (\ref{uph}) is exactly (\ref{upc}).

\end{document}